\documentclass[12pt]{article}
\usepackage{amsmath,enumerate,amsfonts,amssymb,color,graphicx,amsthm}

\usepackage{multirow,caption,subcaption}
\usepackage{hyperref}
\usepackage{todonotes}
\usepackage[normalem]{ulem}

\usepackage{cite}

\numberwithin{equation}{section}

\def\Bbb{  \mathbb }
\newcounter{marnote}

\begin{document}
\newtheorem{lem}{Lemma}
\newtheorem{rem}{Remark}
\newtheorem{question}{Question}
\newtheorem{prop}{Proposition}
\newtheorem{cor}{Corollary}
\newtheorem{thm}{Theorem}[section]
\newtheorem{definition}{Definition}
\newtheorem{openproblem}{Open Problem}
\newtheorem{conjecture}{Conjecture}

\newenvironment{dedication}
  {
   \vspace*{\stretch{.2}}
   \itshape             
  }
  {
   \vspace{\stretch{1}} 
  }

\title
{Symmetry of hypersurfaces and the Hopf Lemma}
\author{
YanYan Li\thanks{Partially
 supported by
        NSF grants DMS-1501004, DMS-2000261, and Simons Fellows  Award 677077.
}
\thanks{Email:  yyli@math.rutgers.edu.}
\\
Department of Mathematics\\
Rutgers University\\
110 Frelinghuysen Road\\
Piscataway, NJ 08854\\
}
\date{ }
\maketitle

\centerline{ \it
Dedicated  to  Neil Trudinger
on his  80th birthday
 with friendship and admiration
}

\begin{abstract}
A classical theorem of A.D. Alexandrov says that a
connected compact smooth hypersurface   in Euclidean space with constant 
mean curvature must be a sphere.
We give exposition to some results
 on symmetry properties of hypersurfaces with ordered mean curvature and associated variations of the Hopf Lemma.
Some open problems will be discussed.
\end{abstract}

\setcounter {section} {-1}

\section{Introduction}

H. Hopf  established in \cite{H}
 that an immersion of a topological $2$-sphere in $\mathbb R^3$
with constant mean curvature must be a standard sphere.
He also made
 the conjecture that the conclusion holds for all immersed connected closed hypersurfaces in $\mathbb R^{n+1}$ with constant mean curvature.
A.D. Alexandrov  proved in \cite{A}  that if
$M$ is an embedded connected closed hypersurface with constant mean
curvature, then $M$
 must be a standard sphere.
If $M$ is immersed instead of embedded,  the conclusion does not hold
in general,
as shown by W.-Y. Hsiang in \cite{Hs}
 for $n\ge 3$ and by Wente in \cite{W}  for $n=2$.
A.
 Ros in \cite{R} gave a different proof for the theorem of Alexandrov
making use of the variational properties
 of the mean curvature.

In this note, we give exposition to some
results in
\cite{Li}-\cite{LYY}. 
It is suggested that the reader 
 read the introductions of \cite{LN1}, \cite{LN2} and \cite{LYY}.

Throughout the paper $M$ is a smooth compact connected embedded hypersurface in $\mathbb R^{n+1}$,  
$k(X)=(k_1(X), \cdots, k_n(X))$ denotes the principal curvatures of $M$ at $X$ with respect to the inner normal, and
the mean curvature of $M$ is
$$
H(X) := \frac 1n\left[ k_1(X)+\cdots + k_n(X)\right].
$$
We use  $G$ to  denote the open bounded set bounded by $M$.

Li
proved in\cite{Li}
the following result saying
  that if the mean curvature $H: M\to \mathbb R$
 has a Lipschitz extension $K:  \mathbb R^{n+1}  \to \mathbb R$ which is 
monotone in the $X_{n+1}$ direction, then
$M$ is symmetric about a hyperplane
$X_{n+1}=c$.

\begin{thm} (\cite{Li})
 Let M be a 
smooth  compact connected embeded
 hypersurface without boundary embedded in
$\mathbb R^{n+1}$,  and let $K$ be a Lipschitz function in $\mathbb R^{n+1}$ 
satisfying
\begin{equation}
K(X', B)\le K(X', A),\quad \forall\ 
X'\in \mathbb R^n, \ A\le B.
\label{1new}
\end{equation}
Suppose that at each point $X$ of $M$ the mean curvature $H(X)$ equals $K(X)$.  Then $M$ is symmetric about a hyperplane 
$X_{n+1}= c.
$
\label{TheoremA}
\end{thm}

In \cite{Li}, $K$ was assumed to be $C^1$ for the above result, but the
 proof there only needs $K$ being Lipschitz.

Li and Nirenberg then considered
in \cite{LN1} and \cite{LN2}
 the more general question in which the condition
$H(X)=K(X)$ with $K$ satisfying (\ref{1new}) is replaced by the weaker,
 more natural,  condition:

\noindent 
{\bf Main Assumption.}
\ For any two points $(X', A), (X', B)  \in M$ satisfying
$A\le B$  
 and that $ \{(X',\theta A+(1 - \theta )B):0 \le \theta  \le 1\}$
lies in $\overline G$, we have
\begin{equation}
H(X', B)\le H(X',A).
\label{main}
\end{equation}

They
 showed in \cite{LN1}
 that this assumption alone is not enough
 to guarentee the symmetry of $M$
about some hyperplane $X_{n+1}=c$.
 The  mean curvature  $H: M\to \mathbb R$ of the
 counterexample
 constructed in [\cite{LN1}, Fig. 4]    has  a monotone extension
$K:  \mathbb R^{n+1}  \to \mathbb R$ which is $C^\alpha $
 for every $0<\alpha<1$, but  fails to be  Lipschitz.
The counterexample actually satisfies (\ref{main}) with an equality. 
They also constructed a  counterexample [\cite{LN1}, Section 6]  
showing that the inequality (\ref{main}) does not imply a pairwise equality.

 A conjecture was made in \cite{LN2}
after the introduction of 

\noindent {\bf 
Condition S.}\ 
$ M $ stays on one side of any hyperplane parallel to 
the $X_{n+1}$  axis that is tangent to $M$.

\medskip

\begin{conjecture}  (\cite{LN2})\
\label{ConjectureA} 
Any smooth compact connected embedded hypersurface $M$ in $\mathbb R^{n+1}$ 
satisfying 
the Main Assumption and Condition S must be symmetric about a hyperplace
$X_{n+1}=c$.
\end{conjecture}

  The conjecture for $n=1$ was proved in \cite{LN1}.
For 
 $n\ge 2$, they introduced
 the following  
condition:

\medskip

\noindent {\bf Condition T.}\
 Every line parallel to the
$X_{n+1}$-axis that is tangent to $M$ has contact of finite order.

\medskip

Note that if  $M$ is real analytic then Condition T is
automatically  satisfied.

They proved in [\cite{LN2}, Theorem 1] that
 $M$ is symmetric about a  hyperplane $X_{n+1}=c$ 
under the Main Assumption,  Condition S and T, and a
local convexity condition
near points where the  tangent planes are parallel to the $X_{n+1}$-axis.
For convex $M$, their 
 result is

\medskip

\begin{thm} (\cite{LN2})
\label{TheoremB}
Let $M$ be a smooth compact convex  hypersurface in $\mathbb R^{n+1}$ satisfying the Main Assumption and Condition  T.
Then $M$ must be symmetric about a  hyperplane $X_{n+1}=c$.
\end{thm}

\medskip

The theorem of Alexandrov is  more general in that one can replace the mean curvature by a wide class of symmetric functions of the principal curvatures.
Similarly,  Theorem \ref{TheoremA}  and 
 Theorem \ref{TheoremB}
(as well as the more general    [\cite{LN2}, Theorem 1]) still hold 
when the mean curvature function is replaced by  
 more general curvature functions.

Consider a triple $(M, \Gamma, g)$: Let $M$
be a compact connected $C^2$
hypersurface without boundary embedded in
$\mathbb R^{n+1}$,  
and let  $g(k_1, \cdots, k_n)$ be a $C^3$
function,
symmetric in $(k_1, \cdots, k_n)$, defined in an open convex
 neighborhood $\Gamma$ of
$\{ (k_1(X), \cdots, k_n(X))\ |\
X\in M\}$, and satisfy 
\begin{equation}
\frac {\partial g}{\partial k_i}(k)>0,\ \ \ 1\le i\le n\ \ \ \
\mbox{and}\ \ \ \ \frac {\partial ^2 g}{ \partial k_i\partial k_j}(k)\eta^i\eta^j
\le 0,\qquad \forall\ k\in \Gamma \ \mbox{and} \ \eta\in \mathbb R^n.
\label{general}
\end{equation}

For convex $M$, their
  result ([\cite{LN2}, Theorem 2]) is
as follows.
\medskip

\begin{thm} (\cite{LN2})\
\label{TheoremC} Let the triple  $(M,\Gamma,g)$ satisfy (\ref{general}). 
In addition, we assume that   $M$ is convex and 
satisfies
Condition T and 
   the Main Assumption  with inequality (\ref{main}) replaced by
\begin{equation}
g(k(X', B))\le g(k(X',A)).
\label{main2}
\end{equation}
   Then $M$ must be symmetric about a hyperplane $X_{n+1}=c$.
\end{thm}

For $1\le m\le n$, let
$$
\sigma_m(k_1, \cdots, k_n)= \sum_{ 1\le i_1<\cdots<i_m\le n} k_{i_1}\cdots k_{i_m}
$$
be the $m$-th elementary symmetric function, and let
$$
g_m:= (\sigma_m)^{\frac 1m}.
$$
It is known that $g=g_m$ satisfies the above properties in
$$
\Gamma_m := \{ (k_1, \cdots, k_n)\in \mathbb R^n\ |\
\sigma_j(k_1, \cdots, k_n)>0\ \mbox{for}\
1\le j\le m\}.
$$
It is known that
$\Gamma_1=\{ k\in \mathbb R^n\ |\ k_1+\cdots + k_n>0\}$,
 $\Gamma_n=\{ k\in \mathbb R^n\ |\ k_1, \cdots, k_n>0\}$,
$\Gamma_{m+1}\subset \Gamma_m$, and
 $\Gamma_m$ is the connected component of
$\{ k\in \mathbb R^n\ |\ \sigma_m(k)>0\}$ containing
$\Gamma_n$.

\medskip

The method of proof of Theorem \ref{TheoremB} and 
\ref{TheoremC}
(as well as
  the more general [\cite{LN2},  Theorem 1 and 2]) 
begins as in that of the theorem of  
 Alexandrov, using the method of moving planes.
Then, as indicated in the introduction of 
\cite{LN1}, one is led to the need for 
variations
 of the classical Hopf Lemma.
The Hopf Lemma is a local result.  The needed variant 
of the Hopf Lemma to prove  Theorem \ref{TheoremB} (and Conjecture
\ref{ConjectureA})   
was raised as an open problem ([\cite{LN2}, Open Problem 2])
which remains open.
The proof of Theorem \ref{TheoremB} and
\ref{TheoremC}
(as well as
  the more general [\cite{LN2},  Theorem 1 and 2]) was 
based on the maximum principle, but  also   used
  a global argument.

In a recent paper \cite{LYY},
Li, Yan and Yao proved
Conjecture \ref{ConjectureA}
  using a method
different from that of \cite{LN1} and \cite{LN2}, exploiting
 the variational properties of 
the mean curvature. 
In fact, they  proved the symmetry result under a 
slightly  weaker assumption than Condition S:

\medskip

\noindent {\bf Condition S'.}\
There exists some constant $r>0$,
 such that for every
$\overline X=(\overline X', \overline X_{n+1})\in M$ with a horizontal
 unit outer normal (denote it by $\bar \nu =(\bar \nu', 0))$,
 the vertical cylinder
$|X'-(\overline X'+r\bar \nu')|=r$ 
has an empty intersection with $G$.
 ($G$  is the bounded open set in
$\mathbb R^{n+1}$  bounded by the hypersurface $M$.)

\begin{thm}  (\cite{LYY})\ 
\label{TheoremD}
 Let $M$
be a compact connected $C^2$
 hypersurface without boundary embedded in
$\mathbb R^{n+1}$,  
 which satisfies both the Main Assumption and Condition
S'.
 Then $M$
must be symmetric about a hyperplane $X_{n+1}=c$.
\end{thm}

Here are two conjectures, in increasing strength.

\begin{conjecture}
For $n\ge 2$ and $2\le m\le n$, let 
$M$
be a compact connected $C^2$
hypersurface without boundary embedded in
$\mathbb R^{n+1}$ satisfying Condition
S (or the slightly weaker Condition S') and
$\{ (k_1(X), \cdots, k_n(X))\ |\
X\in M\}\subset \Gamma_m$.
We assume that  $M$ satisfies the Main Assumption with inequality (\ref{main}) replaced by
\begin{equation}
\sigma_m(k(X', B))\le \sigma_m(k(X',A)).
\label{main3}
\end{equation}
 Then    $M$ must  be symmetric about a  hyperplane $X_{n+1}=c$.
\label{open1}
\end{conjecture}

The next one is for more general curvature functions.

\begin{conjecture}
For $n\ge 2$, let the triple   $(M, \Gamma, g)$ satisfy (\ref{general}).
In addition, we assume that 
$M$
satisfies Condition
S (or the slightly weaker Condition S') and the Main Assumption 
with
 inequality (\ref{main}) replaced by  (\ref{main2}). 
Then    $M$ must  be symmetric about a  hyperplane $X_{n+1}=c$.
\label{open1new}
\end{conjecture}

The above two conjectures  are open even for convex $M$.

 Conjecture \ref{open1} can be approached by two ways. One is by the 
method of moving planes, and this leads to the study of variations of the 
Hopf Lemma.
Such variations of the Hopf Lemma are of
 its own interest.
A number of open problems and conjectures
on such variations of the Hopf Lemma has been discussed in
\cite{LN1}-\cite{LN3}. 
For related works, see \cite{PWY} and  \cite{SB}.
 We will give some discussion on this in Section 1.

Conjecture \ref{open1} can also be approached by using the
variational properties of the higher order mean curvature (i.e. the $\sigma_m$-curvature).
If  the answer to Conjecture \ref{open1}
is affirmative, then the inequality in (\ref{main3}) must be an equality.
This curvature equality
was proved in \cite{LYY2}, using the variational properties
of the $\sigma_m$-curvature:

\begin{lem} (\cite{LYY2})
\label{lemma2new}
For $n\ge 2$ and $2\le m\le n$, let
$M$
be a compact connected $C^2$
hypersurface without boundary embedded in
$\mathbb R^{n+1}$ satisfying Condition
S'.
We assume that $M$ satisfies the Main Assumption, with inequality (\ref{main}) replaced by
(\ref{main3}).
Then (\ref{main3}) must be an equality for every pair of points.
\end{lem}

The proof of Theorem  
\ref{TheoremD} and Lemma \ref{lemma2new} will be sketched in Section 2.

\medskip

We have discussed in the above symmetry properties  of
  hypersurfaces in the Euclidean space.
It is also 
interesting to study
symmetry   properties of hypersurfaces under ordered curvature assumptions 
in the 
hyperbolic space, including
the study of 
 the counter part
 of Theorem \ref{TheoremA}, Theorem
\ref{TheoremD}, and Conjecture
\ref{open1} in the hyperbolic space.
Extensions of the  Alexandrov-Bernstein theorems in the hyperbolic space
were given by Do Carmo and Lawson in \cite{DL};
see also Nelli \cite{N} 
for a survey on  Alexandrov-Bernstein-Hopf theorems.
  
\section{Discussion on Conjecture \ref{open1} and 
 the proof of Theorem \ref{TheoremC} 
 }

Let 

\begin{equation}
\Omega=\{(t,y)\ |\ y\in \Bbb R^{n-1}, |y|<1,
0<t<1\},
\label{D1-1new}
\end{equation}
$$
u, v\in C^\infty(\overline \Omega),
$$
$$
u\ge v\ge 0,\qquad \mbox{in}\ \Omega,
$$
$$
u(0,y)=v(0,y),\quad \forall\ |y|<1;
\qquad u(0,0)=v(0,0)=0,
$$
$$
u_t(0,0)=0,
$$
$$
u_t>0,\qquad \mbox{in}\ \Omega.
$$
We use $k^u(t,y)=(k_1^u(t,y), \cdots, k_n(t,y))$ to 
denote the principal curvatures
 of the graph of $u$ at $(t,y)$.
Similarly, $k^v=(k_1^v, \cdots, k_n^v)$ denotes the  principal 
 curvatures
 of the graph of $v$.

Here are two plausible variations of the Hopf Lemma.
\begin{openproblem}
For $n\ge 2$ and $1\le m\le n$, let $u$ and $v$ satisfy the above.
Assume
$$
\left\{
\begin{array}{l}
\mbox{whenever}\ u(t,y)=v(s,y), 0<s<1, |y|<1,\ \mbox{then there}\\
\sigma_m(k^u)(t,y)\le \sigma_m(k^v)(s,y).
\end{array}
\right.
$$
 
Is it true that either
\begin{equation}
u\equiv v\ \ \mbox{near}\ (0,0)
\label{openproblem}
\end{equation}
or
\begin{equation}
v\equiv 0\ \ \mbox{near}\ (0,0)?
\label{open2}
\end{equation}
\label{OP3}
\end{openproblem}

A weaker version is

\medskip

\begin{openproblem}
 In addition to the assumption
in Open Problem 1, we further assume that
$$
w(t,y):=
\left\{
\begin{array}{ll}
v(t,y),& t\ge 0, |y|<1\\
u(-t, y),&  t<0, |y|<1
\end{array}
\right.
\ \mbox{is}\ C^\infty\ \mbox{in}\
\{(t,y)\ |\ |t|<1, |y|<1\}.
$$
Is it true that either (\ref{openproblem}) or (\ref{open2}) holds?
\label{OP4}
\end{openproblem}

Open Problem \ref{OP3} and \ref{OP4}  for $m=1$ are exactly the same as
[\cite{LN2},  Open Problem 1 and 2],
where it was pointed out that an affirmative 
answer to Open Problem \ref{OP4} for $m=1$ would yield
 a proof of Conjecture \ref{ConjectureA}
 by modification of the arguments in \cite{LN1} and \cite{LN2}.
This applies to $2\le m\le n$ as well: 
An affirmative  
answer to Open Problem \ref{OP4} for some $2\le m\le n$ would yield  a proof of
Conjecture \ref{open1} (with Condition S)
 for the $m$.

As mentioned earlier,
the answer to Open Problem \ref{OP3} for n=1 is yes, and was proved in \cite{LN1}.
For $n\ge 2$, a number of conjectures and open problems
on plausible variations 
 to the Hopf Lemma  were given in 
\cite{LN1}-\cite{LN3}.
The study of such variations of the Hopf Lemma can first be made for the 
Laplace operator instead of the curvature operators.
The following was studied  in \cite{LN3}.

Let $u\ge v$ be in $C^\infty(\overline \Omega)$ where $\Omega$ is given by (\ref{D1-1new}).
Assume that
$$
u>0, \ v>0,\ u_t>0\quad
\mbox{in}\ \Omega
$$
and
$$
u(0, y)=0\quad\mbox{for}\ |y|<1.
$$
We impose a main condition for the Laplace operator:
$$
\mbox{whenever}\ 
u(t,y)=v(s,y)\ \mbox{for}\
0<t\le s<1, \mbox{there}\
\Delta u(t,y)\le \Delta v(s, y).
$$
Under some conditions we wish to conclude that
\begin{equation}
u\equiv v\ \ \mbox{in}\ \Omega.
\label{7}
\end{equation}

The following two conjectures, in decreasing strength, were given
in \cite{LN3}.

\begin{conjecture}
Assume, in addition to the above, that 
\begin{equation}
u_t(0, 0)=0.
\label{8}
\end{equation}
Then (\ref{7}) holds:
$$
u\equiv v\ \ \mbox{in}\ \Omega.
$$
\end{conjecture}

\begin{conjecture}
In addition to (\ref{8}) assume that
$$
u(t,0)\ \mbox{and}\
v(t,0)\ \mbox{vanish at}\
t=0\ \mbox{of finite order}.
$$
Then 
$$
u\equiv v\ \ \mbox{in}\ \Omega.
$$
\end{conjecture}

Partial results  were given in \cite{LN3} concerning these 
conjectures.
On the other hand,  the conjectures  remain largely open.

\section{Discussion on
 Conjecture \ref{open1} and 
 the proof of Theorem \ref{TheoremD} 
}

Theorem \ref{TheoremD} was proved in \cite{LYY} by making use of
the variational properties
 of the mean curvature operator.
We sketch the proof of Therem \ref{TheoremD}  below, see  \cite{LYY}  for
 details.

For any smooth,  closed hypersurface $M$ embeded
in $\mathbb R^{n+1}$, let
$V:
\mathbb R^{n+1}\to \mathbb R^{n+1}$ be a  smooth vector field.
Consider, for $|t|<1$,
\begin{equation}
M(t):= \{ x+tV(x) \ |\ x\in M\},
\label{M}
\end{equation}
and
$$
S(t):= \int_{ M(t) }d\sigma= \mbox{area of}\ M(t).
$$
It is well known that
\begin{equation}
\frac{d}{dt}S(t)\bigg|_{ t=0}=-\int_M
V(x)\cdot \nu(x) H(x) d\sigma(x),
\label{var}
\end{equation}
where $H(x)$ is the mean curvature of $M$ at $x$ with respect to the inner
unit normal  $\nu$.

  Define the projection map   $\pi: (x', x_{n+1}) \to x'$, and set
 $R:=\pi (M)$.

 Condition S' assures that 
 $\nu(\bar x)$, $\bar x\in M$,  is horizontal
 iff $\bar x'\in \partial R$;
  $\partial R$  is $C^{1,1}$
(with $C^{1,1}$ normal  under control); and
 $$
M=M_1\cup M_2\cup \widehat M,$$
where
$M_1, M_2$ are respectively graphs of functions   $f_1, f_2: R^\circ\to R,$
$f_1, f_2\in C^2(R^\circ), f_1>f_2$ in
$R^\circ$, and
$\widehat  M:=\{ (x', x_{n+1})\in M\ |\ x'\in \partial R\}\equiv M\cap \pi^{-1}(\partial R).$
Note that $f_1, f_2$ are not in $C^0(R)$ in general.

\begin{lem}
\begin{equation}
H(x', f_1(x'))= H(x', f_2(x'))\quad  \forall\ x'\in R^\circ.
\label{equality}
\end{equation}
\label{lemma1a}
\end{lem}

\noindent {\bf Proof.}\ 
Take $V(x)=e_{n+1}=(0, ..., 0, 1)$,
and let $M(t)$ and
$
S(t)$
be defined as above with this choice of $V(x)$.
Clearly, $S(t)$ is  independent of $t$.
So
we have, using (\ref{var}) and the order assumption on the mean curvature,
that
\begin{eqnarray}
0&=&\frac{d}{dt}S(t)\bigg|_{ t=0}=-\sum_{i=1}^2
\int_{M_i}
e_{n+1} \cdot \nu(x) H(x) d\sigma(x)\nonumber\\
&=& -\int_{R^\circ}
\left[H(x', f_1(x'))-H(x', f_2(x'))\right]dx'\ge 0.
\label{3.5}
\end{eqnarray}
Using  again the order assumption on the mean curvature we obtain
The curvature equality (\ref{equality}).

\medskip

 For any $v\in
C^\infty(R^{n})$, let  
$V(x):=v(x')e_{n+1}$, and let $M(t)$ and
$
S(t)$
be defined as above with this choice of $V(x)$.
We have, using (\ref{var}) and  (\ref{equality}), that
\begin{eqnarray}
0&=&\frac{d}{dt}S(t)\bigg|_{ t=0}=-\sum_{i=1}^2
\int_{M_i}
v(x') e_{n+1} \cdot \nu(x) H(x) d\sigma(x)\nonumber\\
&=& -\int_{R^\circ}
v(x') \left[H(x', f_1(x'))-H(x', f_2(x'))\right]dx'= 0.
\label{One}
\end{eqnarray}

Theorem \ref{TheoremD} is proved by contradition as follows:
If $M$ is not symmetric about a hyperplane, then
  $\nabla (f_1+f_2)$ is not identically zero.  We will  find a
particular $V(x)=v(x') e_{n+1}$,  $v\in C^2_{loc}(R^\circ)
$, 
 to make
$$
\frac{d}{dt}S(t)\bigg|_{ t=0}\ne 0,
$$
which contradicts to (\ref{One}).
 
 Write
$$
S(t)= \sum_{i=1}^2
\int_{ R^\circ}
\sqrt{
1+ |\nabla [ f_i(x') +v(x') ]t|^2 } dx'+\widehat S,
$$
where $\widehat S$, the area of the vertical part of $M$,
  is independent of $t$ (since $v$ is zero near 
$\partial R$, so the  vertical part of $M$ is not moved).

\medskip

A calculation gives
$$
\frac{d}{dt}S(t)\bigg|_{ t=0}
=
\int_{ R^\circ}
 \sum_{i=1}^2
\left[\nabla A(\nabla f_1(x'))- \nabla A(-\nabla f_2(x'))\right]
\cdot \nabla v(x')dx',
$$
where
$$
A(q):=\sqrt{ 1+|q|^2},\ \ q\in R^n.
$$

We know
that $$
\nabla A(q)= \frac q{  \sqrt{ 1+|q|^2} }\quad \mbox{and} 
\quad \nabla^2 A(q)\ge (1+|q|^2)^{ -3/2} I>0\ \ \forall q.
$$
So 
 $[\nabla A(q_1)-\nabla A(q_2)] \cdot (q_1-q_2)>0$
for any $q_1\ne q_2$.

\medskip

If $\nabla (f_1+f_2)\in C^2_{loc}(R^\circ)u\setminus \{0\}$, we would take
 $v=\nabla (f_1+f_2)$ and obtain
\begin{eqnarray*}
\frac{d}{dt}S(t)\bigg|_{ t=0}
&=&
\int_{ R^\circ}
\left[\nabla A(\nabla f_1(x'))- \nabla A(-\nabla f_2(x'))\right]
\cdot \nabla v(x')dx'\\
&=&
\int_{ R^\circ}
\left[\nabla A(\nabla f_1(x'))- \nabla A(-\nabla f_2(x'))\right]
\cdot  [\nabla f_1(x')+ \nabla f_2(x')]dx'\\
&>&
0.
\end{eqnarray*}

In general,  $\nabla (f_1+f_2)$ would not be in $C^2_{loc}(R^\circ)$.
It turns out that   Condition S' allows us to do a smooth
 cutoff near
$\partial R$, and conclude the proof.
We skip the crucial details, which can be found in the last few pages of
  \cite{LYY}.

\medskip

\bigskip

Now we give the

\noindent{\bf Proof of Lemma \ref{lemma2new}.}\ The proof is similar to that of Lemma \ref{lemma1a}, see also the proof of [\cite{LYY}, Proposition 3] for more details.
We still   take $V(X)=e_{n+1}$ 
and let $M(t)$ be as in (\ref{M}).  Consider
$$
S_{m-1}(t):=
\int_{ M(t) }
\sigma_{m-1} (x) d\sigma.
$$
Clearly, $S_{m-1}(t)$ is independent of $t$.

The variational properties of higher order curvature
[\cite{Reilly}, Theorem B] gives
\begin{eqnarray*}
\frac{d}{dt}S_{m-1}(t)\bigg|_{ t=0}=-m\int_M
V(x)\cdot \nu(x) \sigma_m(x) d\sigma(x),
\end{eqnarray*}
thus the same argument as (\ref{3.5}) yields 
\begin{eqnarray*}
0&=&\frac{d}{dt}S_{m-1}(t)\bigg|_{ t=0}=-m\int_M
V(x)\cdot \nu(x) \sigma_m(x) d\sigma(x)\\
&=& -\int_{R^\circ}
\left[\sigma_m(x', f_1(x'))-\sigma_m(x', f_2(x'))\right]dx'\ge 0.
\end{eqnarray*}
We deduce from the above, using the curvature inequality  (\ref{main3}),
that the equality in  (\ref{main3}) must hold for every pair of points.
Lemma \ref{lemma2new} is proved.

\end{document}